\documentclass[11pt]{amsart}

\usepackage{amsmath,amsthm,amscd,euscript}
\def \z{\mathbb Z}

\def \f{\mathcal F}

\DeclareMathOperator{\MOD}{mod} 
 \DeclareMathOperator{\Pol}{Pol}
\DeclareMathOperator{\Exp}{Exp}\DeclareMathOperator{\ord}{ord}

\newtheorem{theorem}{Theorem}[section]

\newtheorem{proposition}[theorem]{Proposition}

\theoremstyle{remark}
\newtheorem{remark}[theorem]{Remark}
\theoremstyle{definition}
\newtheorem{definition}[theorem]{Definition}

\newtheorem{problem}{Problem}

\title{On examples of difference operators for $\{0,1\}$-valued functions
over finite sets}
\author{Oleg~Karpenkov}
\date{31  October 2006}
\thanks{Partially supported by NWO-RFBR 047.011.2004.026
(RFBR 05-02-89000-NWO\_a) grant, by
RFBR SS-1972.2003.1 grant, by RFBR 05-01-02805-CNRSL\_a grant,
and by RFBR grant 05-01-01012a.}


\email[Oleg Karpenkov]{karpenk@mccme.ru}

\address{Mathematisch Instituut, Universiteit Leiden,
P.O. Box 9512, 2300 RA Leiden, The Netherlands}

\begin{document}
\input epsf

\maketitle


\section{Introduction and basic definitions}

Recently V.~I.~Arnold have formulated a geometrical concept of
monads and apply it to the study of difference operators on the
sets of $\{0,1\}$-valued sequences of length $n$.
In~\cite{ArnC1}--\cite{ArnC4} he made first steps in the study of
this subject and formulated many nice questions. In~\cite{Gar}
A.~Garber showed an algorithm that gives a description of the
combinatorial structure of monads for difference operators and
answered many questions of V.~I.~Arnold. In the present note we
show particular examples of these monads and indicate one
question arising here.

The author is grateful to V.~I.~Arnold for useful remarks and
discussions and Mathematisch Instituut of Universiteit Leiden for
the hospitality and excellent working conditions.

\vspace{2mm}

A {\it monad} by V.~I.~Arnold is a map of a finite set into
itself. Suppose $M:S\to S$ is an arbitrary monad. It is naturally
to associate an oriented graph to the monad $M$. Its vertices
coincide with elements of $S$, and the set of its edges is the
set of ordered pairs $(x,M(x))$. We denote such graph by $G(M)$
The idea of V.~I.~Arnold is to study the combinatorial geometry of
graphs for monads. He proposed to start with one important
example.

Consider any positive integer $n$, and take the set
$A_n=\{1,\ldots, n\}$. Denote by $\f_2(A_n)$ the vector space of
$\z_2$-valued functions on $A_n$. Consider a ``differential'' {\it
difference} operation $\Delta$, defined as follows:
$$
(\Delta f)(x)= \left\{
\begin{array}{ll}
f(x)+f(x+1), & \hbox{if $x\ne n$}\\
f(n)+f(1), & \hbox{if $x=n$}\\
\end{array}
. \right.
$$

\begin{figure}[h]
$$\epsfbox{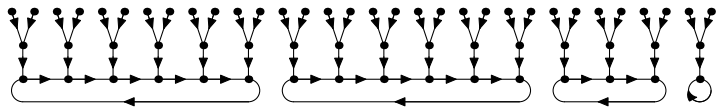}$$
\caption{An example of a monad $G(\Delta)$ for the case of
$A_{6}$.}\label{monada.1}
\end{figure}

On Figure~\ref{monada.1} we show an example of a monad for the
difference operator for the set $A_6$.

\section{A few words about Arnold's complexity}

Let us briefly describe the concept of Arnold's functional
complexity. Further we will need the additivity property of
$\Delta$:
$$
\Delta(f+g)=\Delta(f)+\Delta(g).
$$

The simplest functions on the set $A_n$ are polynomials. The set
of all solutions of the functional equation $\Delta ^{k} (f)=0$
is called the {\it set of polynomials of degree less than $k$}, we
denote it by $\Pol_{k-1}$. Suppose $f\in \Pol_{k}$ and $f\notin
\Pol_{l}$ for $l=0,\ldots, k{-}1$, then $f$ is called a {\it
polynomial of degree $k$}. Denote by $\Pol (A_n)$ the set of all
polynomials on $A_n$.

Actually, the set $\Pol (A_n)$ is a vector space. If $n=2^l m$
where $m$ is odd, then the space $\Pol (A_n)$ is
$2^l$-dimensional and contains $2^{2^l}$ elements. In particular,
if $n=2^l$, then $\f_2(A_n)=\Pol (A_n)$.

Now we define another set of nice functions. Consider the
functional equation $\Delta ^{k} (f)=f$. The set of all solutions
of such equation is called the {\it set of special rational
exponential polynomials of orders that divide $k$}, we denote it
by $\Exp_{k}$. Let $\Exp_{0}=\{0\}$. Suppose $f\in \Exp_{k}$ and
$f\notin \Exp_{l}$ for $l=0,\ldots, k{-}1$, then $f$ is called a
{\it special rational exponential polynomial of order $k$}, or
{\it exp-polynomial} for short. Denote by $\Exp (A_n)$ the set of
all exp-polynomials on $A_n$. Note that the set $\Exp (A_n)$ is a
vector space.

\begin{proposition}
Any function $f$ of $\f_2(A_n)$ can be uniquely written in the
form $f=p+r$, where $p$ is a polynomial and $r$ is an
exp-polynomial, or in other worlds
$$
\f_2(A_n)=\Pol (A_n)\oplus \Exp (A_n).
$$
\end{proposition}

\begin{definition}
Consider an arbitrary function of $\f_2(A_n)$. Let $f=p+r$, where
$p$ is a polynomial and $r$ is an exp-polynomial. We say that {\it
degree} of $f$ is degree of $p$ and denote it by $\deg(f)$. We
say that {\it order} of $f$ is order of $r$ and denote it by
$\ord(f)$.
\end{definition}

V.~I.~Arnold proposed the following definition the notion of
functional complexity.
\begin{definition}
A function $f_1$ is said to be {\it more complicated $($in the
sense of Arnold$)$} than $f_2$ if either $\ord(f_1)>\ord(f_2)$, or
$\ord(f_1)=\ord(f_2)$ and $\deg(f_1)>\deg(f_2)$.\\
If $\ord(f_1)=\ord(f_2)$ and $\deg(f_1)=\deg(f_2)$ then the
functions $f_1$ and $f_2$ are said to be of the {\it same
complexity $($in the sense of Arnold$)$}.
\end{definition}

\begin{remark}
All the above easily can be generalized to the case of
$\z_N$-valued functions for an arbitrary positive integer $N>2$.
\end{remark}

\section{Some examples of $G(\Delta)$}

\subsection{A list of examples.} Suppose $n=2^lm$. The set of polynomials is a $2$-valent
symmetric tree of $r=2^{2^l}$ elements. Denote this tree by
$T_r$. Denote also a cycle of $s$ elements by $O_s$.

Each connected component of the graph $G(\Delta)$ contains a
cycle. Denote its length by $s$. To each vertex of the cycle it is
attached a tree equivalent to $T_r$ (as on Fig.~\ref{monada.1} for
the case $n=6$). Denote such component by $O_s{*}T_r$.

Let us enumerate connected components for the graphs of the sets
$\f_2(A_n)$ where $n\le 25$. Expression $k(O_s*T_r)$ means that
there are $k$ components of the type $O_s*T_r$. The graph
$G(\Delta)$ for the set $\f_2(A_n)$ contains\\
in the case of $n=1$: $O_1{*}T_{2}$;\\
in the case of $n=2$: $O_1{*}T_{4}$;\\
in the case of $n=3$: $O_1{*}T_{2}$, $O_3{*}T_{2}$;\\
in the case of $n=4$: $O_1{*}T_{16}$;\\
in the case of $n=5$: $O_1{*}T_{2}$, $O_{15}{*}T_{2}$;\\
in the case of $n=6$: $O_1{*}T_{4}$, $O_3{*}T_{4}$, $2(O_6{*}T_{4})$;\\
in the case of $n=7$: $O_1{*}T_{2}$, $9(O_7{*}T_{2})$;\\
in the case of $n=8$: $O_1{*}T_{256}$;\\
in the case of $n=9$: $O_1{*}T_{2}$, $O_3{*}T_{2}$, $4(O_{63}{*}T_{2})$;\\
in the case of $n=10$: $O_1{*}T_{4}$, $O_{15}{*}T_{4}$, $8(O_{30}{*}T_{4})$;\\
in the case of $n=11$: $O_1{*}T_{2}$, $3(O_{341}{*}T_{2})$;\\
in the case of $n=12$: $O_1{*}T_{16}$, $O_3{*}T_{16}$, $2(O_6{*}T_{16})$, $20(O_{12}{*}T_{16})$;\\
in the case of $n=13$: $O_1{*}T_{2}$, $5(O_{819}{*}T_{2})$;\\
in the case of $n=14$: $O_1{*}T_{4}$, $9(O_7{*}T_{4})$, $288(O_{14}{*}T_{4})$;\\
in the case of $n=15$: $O_1{*}T_2$, $O_3{*}T_2$, $30(O_5{*}T_2)$, $1082(O_{15}{*}T_2)$;\\
in the case of $n=16$: $O_1{*}T_{2^{16}}$;\\
in the case of $n=17$: $O_1{*}T_2$, $51(O_{85}{*}T_2)$, $240(O_{255}{*}T_2)$;\\
in the case of $n=18$: $O_1{*}T_4$, $O_{3}{*}T_4$, $2(O_{6}{*}T_4)$, $4(O_{63}{*}T_4)$, $518(O_{126}{*}T_4)$;\\
in the case of $n=19$: $O_1{*}T_2$, $27(O_{9709}{*}T_2)$;\\
in the case of $n=20$: $O_1{*}T_{16}$, $O_{15}{*}T_{16}$, $8(O_{30}{*}T_{16})$, $1088(O_{60}{*}T_{16})$;\\
in the case of $n=21$: $O_1{*}T_2$, $O_{3}{*}T_2$, $9(O_{7}{*}T_2)$, $9(O_{21}{*}T_2)$, $16640(O_{63}{*}T_2)$;\\
in the case of $n=22$: $O_1{*}T_4$, $3(O_{341}{*}T_4)$, $1536(O_{682}{*}T_4)$;\\
in the case of $n=23$: $O_1{*}T_2$, $2049(O_{2047}{*}T_2)$;\\
in the case of $n=24$: $O_1{*}T_{256}$, $O_{3}{*}T_{256}$, $2(O_{6}*T_{256})$, $20(O_{12}{*}T_{256})$, $2720(O_{24}{*}T_{256})$;\\
in the case of $n=25$: $O_1{*}T_{2}$, $O_{15}{*}T_{2}$,
$656(O_{25575}{*}T_{2})$.

Denote by $s(n)$ the order of the maximal possible length of
cycles for the $n$-elements sequences. Actually the listed
examples gives the negative answer to the following question of
V.~I.~Arnold: {\it is it true that $(s(n)/n){+}1$ is some power of
$2$}? It is not true, for example, for $n=23$ where $s(23)=2047$.
Here $s(23)$ is $2^{11}{-}1$ itself.

Denote by $]n[$ the set of connected components of graphs
$G(\Delta)$, corresponding to the set $\f_2(A_n)$. The
work~\cite{Gar} of A.~Garber immediately implies the following
identities:
$$
\begin{array}{l}
\hbox{a) } ]3\cdot 2^m[=\left\{O_1{*}T_{2^{2^m}},
O_3{*}T_{2^{2^m}}\right\} \bigcup
\left\{\frac{2^{3\cdot2^k}-2^{4\cdot2^{k-1}}}{3\cdot2^k
\cdot2^{2^k}} (O_{3\cdot 2^k}{*}T_{2^{2^m}}) \Big|k=2,\ldots,
m\right\};\\
\hbox{b) } ]5\cdot 2^m[=\left\{O_1{*}T_{2^{2^m}}, O_{3\cdot
5}{*}T_{2^{2^m}}\right\} \bigcup
\left\{\frac{1}{3}\cdot\frac{2^{5\cdot2^k}-2^{6\cdot2^{k-1}}}{5\cdot2^k
\cdot2^{2^k}} (O_{3\cdot 5\cdot 2^k}{*}T_{2^{2^m}})
\Big|k=2,\ldots,
m\right\};\\
\hbox{c) } ]7\cdot 2^m[=\left\{O_1{*}T_{2^{2^m}},
9(O_7{*}T_{2^{2^m}})\right\} \bigcup
\left\{\frac{2^{7\cdot2^k}-2^{8\cdot2^{k-1}}}{7\cdot2^k
\cdot2^{2^k}} (O_{7\cdot 2^k}{*}T_{2^{2^m}}) \Big|k=2,\ldots,
m\right\};\\
\hbox{d) } ]11\cdot 2^m[=\left\{O_1{*}T_{2^{2^m}}, 3(O_{31\cdot
11}{*}T_{2^{2^m}})\right\} \bigcup \left\{\frac{1}{31}
\frac{2^{11\cdot2^k}-2^{12\cdot2^{k-1}}}{11\cdot2^k \cdot2^{2^k}}
(O_{31\cdot 11\cdot 2^k}{*}T_{2^{2^m}}) \Big|k=2,\ldots, m\right\}\\
\cdots\\
\end{array}
$$

\subsection{Particular case of $\delta$-functions.}
Let us now study the structure of the piece-wise connected
components of the graph $G$ containing so-called {\it
$\delta$-function}. Denote by $\delta_k$ the following function
of $\f_2(A_n)$:
$$
\delta_k (x)= \left\{
\begin{array}{ll}
0, & \hbox{if $x\ne k$}\\
1, & \hbox{if $x=k$}\\
\end{array}
. \right.
$$

In~\cite{Gar} A.~Garber showed that the order of $\delta_k$
coincide with $s(n)$. So the piece-wise connected component of the
graph $G$ containing $\delta_k$ is $O_{s(k)}{*}T_{2^{2^l}}$, and
it does not depend on the choice of $k$. We now write down the
values of $s(n)$ for $n\le 50$ in the following list.

\begin{center}
\begin{tabular}{|c||c|c|c|c|c|c|c|c|c|c|c|c|c|c|c|c|c|c|c|c|}
\hline
n&1&2&3&4&5&6&7&8&9&10&11&12&13&14&15&16&17&18&19&20\\
\hline
$s(n)$ &1&1&3&1&15&6&7&1&63&30&341&12&819&14&15&1&255&126&9709&60\\
\hline
\end{tabular}

\vspace{2mm}

\begin{tabular}{|c||c|c|c|c|c|c|c|c|c|c|c|c|c|c|c|}
\hline
n&21&22&23&24&25&26&27&28&29&30&31&32&33&34&35\\
\hline
$s(n)$&63&682&2047&24&25575&1638&13797&28&475107&30&31&1&1023&510&4095\\
\hline
\end{tabular}

\vspace{2mm}

\begin{tabular}{|c||c|c|c|c|c|c|c|c|c|c|c|c|}
\hline
n&36&37&38&39&40&41&42&43&44&45&46&47\\
\hline
$s(n)$&252&3233097&19418&4095&120&41943&126&5461&1364&4095&4095&8388607\\
\hline
\end{tabular}

\vspace{2mm}

\begin{tabular}{|c||c|c|c|}
\hline
n&48&49&50\\
\hline
$s(n)$ &48&2097151&51150\\
\hline
\end{tabular}

\end{center}

There is a regularity in this sequence for primes $n>2$. We use
the following notation. Denote by $\gamma_2(n)$ the minimal
solution $t$ of the equation $2^t\equiv 1 (\MOD n)$. Then

\begin{center}

\begin{tabular}{|c||c|c|c|c|c|c|c|}
\hline
n&3&5&7&11&13&17&19\\
\hline\hline
$\gamma_2(n)$&2&4&3&10&12&8&18\\
\hline
$s(n)$
&$3(2^{\frac{2}{2}}{-}1)$&$5(2^{\frac{4}{2}}{-}1)$&$2^{3}{-}1$&$11(2^{\frac{10}{2}}{-}1)$&
$13(2^{\frac{12}{2}}{-}1)$&$17(2^{\frac{8}{2}}{-}1)$&$19(2^{\frac{18}{2}}-1)$\\
\hline
\end{tabular}

\vspace{2mm}

\begin{tabular}{|c||c|c|c|c|c|c|c|}
\hline
n&23&29&31&37&41&43&47\\
\hline\hline
$\gamma_2(n)$&11&28&5&36&20&14&23\\
\hline
$s(n)$ &$2^{11}{-}1$&$29(2^{\frac{28}{2}}{-}1)$&$2^{5}{-}1$&
$37(2^{\frac{36}{2}}{-}1)/3$&$41(2^{\frac{20}{2}}{-}1)$&$43(2^{\frac{14}{2}}{-}1)$&$2^{23}{-}1$\\
\hline
\end{tabular}

\end{center}

Denote by $q(n)$ the following function
$$
q(n):=\left\{
\begin{array}{ll}
n(2^{\frac{\gamma_2(n)}{2}}-1), & \hbox{if $\gamma_2(n)$ is even}\\
2^{\gamma_2(n)}-1, & \hbox{if $\gamma_2(n)$ is odd}\\
\end{array}
.
\right.
$$

Note that for all observed primes (except 37) we have $s(n)=q(n)$.

\begin{problem}
Study the behaviour of the maximal length of the cycle. How often
does it coincide with $q(n)$? Is it true that $q(n)$ is always
divisible by $s(n)$?
\end{problem}

\end{document}